\documentclass{amsart}

\usepackage{graphicx, psfrag}

\newtheorem{theorem}{Theorem}[section]
\newtheorem{lemma}[theorem]{Lemma}

\theoremstyle{definition}

\theoremstyle{remark}

\numberwithin{equation}{section}

\begin{document}

\title{A note on circle patterns on surfaces}

\author{Ren Guo}

\address{Department of Mathematics, Rutgers University, Piscataway, NJ, 08854}

\email{renguo@math.rutgers.edu}

\thanks{This work is partially supported by NSF Grant \#0625935.}

\subjclass[2000]{52C26}

\keywords{circle pattern, continuity method, variational
principle.}

\begin{abstract}
In this paper we give two different proofs of Bobenko and
Springborn's theorem of circle pattern: there exists a hyperbolic
(or Euclidean) circle pattern with proscribed intersection angles
and cone angles on a cellular decomposed surface up to isometry
(or similarity).
\end{abstract}

\maketitle

\section {Introduction}

In the work of Bobenko and Springborn \cite{bs}, they produced a
variational principle for circle patterns on surfaces in Euclidean
and hyperbolic geometry. The goal of this paper is to give two
different proofs of Bobenko and Springborn's theorem. First we
will give a direct proof of Bobenko and Springborn's theorem in
Euclidean and hyperbolic geometry by using the continuity method
used by Thurston \cite{th}, Marden and Rodin \cite{mr}. Then we
will use another variational principle following Rivin's strategy
\cite{r1} \cite{r2} to prove Bobenko and Springborn's theorem in
hyperbolic geometry.

Bobenko and Springborn's proof consists of two steps. In the first
step they use the feasible flow to find the necessary and
sufficient condition for the existence of an angle structure
(i.e.coherent angle system). In the second step they introduced an
energy function whose variables are radii of circles and they can
show that the unique critical point of the energy function
corresponds to circle pattern with prescribed intersection angles
and cone angles.

This variational principle is first used by Colin de Verdi\'ere
\cite{cdv} in the proof of Andreev-Thurston's circle packing
theorem in Euclidean and hyperbolic geometry. Br\"agger \cite{b}
found a new energy function for Euclidean circle packing whose
variables are certain angles. This energy function turns out to be
the Legendre transformation of Colin de Verdi\'ere's.

To character the cone metric on a triangulated surface, Rivin
produced a variational principle. Rivin \cite{r2} used the duality
theorem in linear programming to find the necessary and sufficient
condition for the existence of an angle structure. Rivin \cite{r1}
introduced an energy function whose variables are certain angles
and he can show the unique maximal point of the energy function
corresponds to the angle structure induced by a cone metric.
Following Rivin's method, Leibon \cite{le} characterized the
hyperbolic cone metric. For the spherical cone metric, the angle
structure part was solved by the author \cite {g} and the energy
function part was solved by Luo \cite {l1}. We will show that
Rivin's strategy also works in the circle pattern problem.

Our energy function is the Legendre transformation of the energy
function in Bobenko and Springborn \cite{bs}. Springborn \cite{s1}
already found this Legendre transformation. Here we will show this
function can be naturally derived from Luo's principle of ``the
derivative cosine law". Luo \cite{l1} established this principle
and was able to recover the energy function of Colin de
Verdi\'ere, Br\"agger, Rivin, Leibon. See the related work of Luo
\cite{l2} \cite{l3} \cite{l4}.

Springborn \cite{s2} used Rivin's strategy to solve another circle
pattern problem.

To state Bobenko and Springborn's theorem, let's recall some
definitions first. Suppose $\Sigma$ is a closed surface with a
cell decomposition. Let $E, F$ be the sets of edges and faces of
the cell decomposition respectively. If $e\in E, f\in F$, we use
$e<f$ to denote that $e$ is an edge of $f$. Suppose we are given a
\it radius function \rm $r:F\to(0,\infty)$ and an \it intersection
angle function \rm $D:E\to (0,\pi)$. We will use these data to
construct a hyperbolic (or Euclidean) structure on $\Sigma$ with
cone points.

\begin{figure}[htbp]
\begin{center}
\includegraphics[scale=.7]{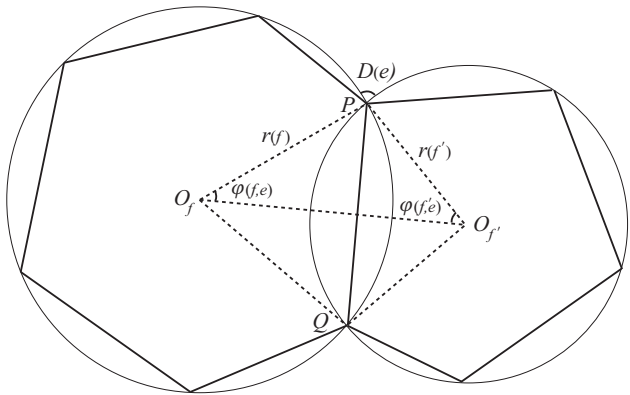}
\end{center}
\caption{\label{Fig:circle}Cell decomposition, quadrilateral,
circle pattern}
\end{figure}

First, choose a point $O_f$ in the interior of each $f.$ As
in Figure \ref{Fig:circle}, if the face $f$ and $f'$ share an edge $e$ with end points
$P$ and $Q,$ we realize the quadrilateral $O_fPO_{f'}Q$ as a hyperbolic (or Euclidean)
quadrilateral with edge lengths $|O_fP|=|O_fQ|=r(f),|O_{f'}P|=|O_{f'}Q|=r(f')$ and
angles $\angle O_fPO_{f'}=\angle O_fQO_{f'}=\pi-D(e).$ For each edge $e$ we construct such a
quadrilateral. By gluing together all the quadrilaterals, we obtain a hyperbolic (or Euclidean)
structure on $\Sigma$ with cone points. Furthermore for each face $f$ there is a circumcircle with
center $O_f$ and radius $r(f).$ And the intersection angle of the circumcircles of face $f$ and $f'$ is $D(e).$
We call this a hyperbolic (or Euclidean) \it circle pattern \rm. In the hyperbolic (or Euclidean) triangle
$O_fPO_{f'}$, denote the angle at $O_f$ by $\varphi(f,e).$ Then the cone angle at point $O_f$ is
$$\sum_{e<f}2\varphi(f,e),$$ where the sum runs over all edges of the face
$f.$ Thus we obtain a \it cone angle function \rm $\Phi: F \to
(0,\infty)$.

On the other hand, if we are given an intersection angle function
$D:E\to (0,\pi)$ and a cone angle function $\Phi: F \to
(0,\infty)$, we try to solve the radius function $r:F \to
(0,\infty)$, equivalently to obtain a circle pattern. Bobenko and
Springborn \cite{bs} obtained the following result by using a
variational principle. Their result can cover the case of surfaces
with boundary. But for simplicity we only consider the closed
surfaces in this paper.

For any subset $X\subseteq F,$ let $E(X)=\{e\in E|e<f\in X\}$ be the set
of edges contained in at least one of the faces in $X.$

\begin{theorem}[Bobenko-Springborn \cite{bs} Theorem 3]\label{e} Given a cellular decomposed surface
$\Sigma$, a function $D:E\to (0,\pi)$ and a function $\Phi: F \to
(0,\infty)$, there exists a unique Euclidean circle pattern with
the the intersection angle function $D$ and the cone angle
function $\Phi$ if and only if

\begin{equation}
\label{ee} \sum_{f\in F}\frac12\Phi(f)= \sum _{e\in E} D(e),
\end{equation}

 and
\begin{equation}
\label{eee} \sum_{f\in X}\frac12\Phi(f)< \sum _{e\in E(X)} D(e)
\end{equation}

for any nonempty subset $X \subset F, X\neq F.$
\end{theorem}

\begin{theorem}[Bobenko-Springborn \cite{bs} Theorem 3]\label{h} Given a cellular decomposed surface
$\Sigma$, a function $D:E\to (0,\pi)$ and a function $\Phi: F \to
(0,\infty)$, there exists a unique hyperbolic circle pattern with
the the intersection angle function $D$ and the cone angle
function $\Phi$ if and only if

\begin{equation}
\label{hh} \sum_{f\in X}\frac12\Phi(f)< \sum _{e\in E(X)} D(e)
\end{equation}

for any nonempty subset $X \subseteq F.$
\end{theorem}

In section 2 we will prove Theorem \ref{e} by using the continuity
method. In section 3 we will prove Theorem \ref{h} by using the
continuity method. In section 4 we will prove Theorem \ref{h} by
using a variational principle.

\section{Continuity method for Euclidean circle pattern}

\begin{proof}[Proof of Theorem \ref{e}] Assume the set of faces of the cell
decomposition is $F=\{f_1,...,f_n\}$. Due to Euclidean similarity,
it is enough to consider the space of radii
$\Delta=\{(r(f_1),...,r(f_n))|r(f_1)+...+r(f_n)=1\}$ which is a
simplex of dimension $n-1.$ For a fixed intersection angle
function $D:E\to(0,\pi)$, if the radius function
$r:F\to(0,\infty)$ is given, we can find the cone angle function
$\Phi:F\to(0,\infty)$. Therefore there is a mapping $P:\Delta\to
\mathbb{R}^n$ sending the radii $(r(f_1),...,r(f_n))$ to its cone
angles $\Phi=(\Phi(f_1),...,\Phi(f_n)).$

Let $\delta=\{\Phi\in\mathbb{R}^n_{>0}|\Phi\ \text{satisfies the
conditions}\ (\ref{ee}),(\ref{eee})\}$ be an $n-1$ dimensional
polytope in $\mathbb{R}^n$. We claim that
$P(\Delta)\subset\delta.$ We need to check that, if $\Phi\in
P(\Delta),$ $\Phi$ satisfies the conditions (\ref{ee}),
(\ref{eee}).

In fact, as in Figure \ref{Fig:circle}, each
Euclidean triangle has inner angles $(\varphi(f,e),\varphi(f',e),\pi-D(e)).$ Then
\begin{equation*}
\begin{split}
\sum_{f\in F}\frac12\Phi(f)
&=\sum_{f\in F}\sum_{e<f}\varphi(f,e)\\
&=\sum_{e\in E,e<f,e<f'}(\varphi(f,e)+\varphi(f',e))\\
&=\sum_{e\in E}D(e).
\end{split}
\end{equation*}
Given a nonempty proper subset of $F$, i.e. $X\subset F, X\neq F,$
the set $E(X)=\{e\in E|e<f\in X\}$ can be decomposed as a union of
two disjoint subsets $E(X)_1\cup E(X)_2,$ where $E(X)_1=\{e\in
E|e<f\in X, e<f'\notin X\}$ and $E(X)_2=\{e\in E|e<f\in X, e<f'\in
X\}$. An edge in $E(X)_1$ is contained in exact one of the faces
in $X$ while an edge in $E(X)_2$ is contained in two of the faces
in $X$. Then
\begin{equation*}
\begin{split}
\sum_{f\in X}\frac12\Phi(f)
&=\sum_{f\in X}\sum_{e<f}\varphi(f,e)\\
&=\sum_{e\in E(X)_1,f\in X}\varphi(f,e)+\sum_{e\in
E(X)_2}(\varphi(f,e)+\varphi(f',e))\\
&=\sum_{e\in E(X)_1,f\in X}\varphi(f,e)+\sum_{e\in E(X)_2} D(e)\\
&<\sum_{e\in E(X)_1}D(e)+\sum_{e\in E(X)_2} D(e)\\
&=\sum_{e\in E(X)}D(e).
\end{split}
\end{equation*}
Next we claim the mapping $P:\Delta\to\delta$ is one to one. In
fact, let $r=(r(f_1),...,r(f_n))$ and $
\overline{r}=(\overline{r}(f_1),...,\overline{r}(f_n))$ be two
distinct points in $\Delta.$ Let $F_0=\{f_i\in
F|r(f_i)>\overline{r}(f_i)\}$. It is a nonempty proper subset. By
the calculation above for a general subset $X$, we have
$$\sum_{f\in F_0}\frac12\Phi(f)=\sum_{e\in E(F_0)_1,f\in F_0}\varphi(f,e)+\sum_{e\in E(F_0)_2} D(e).$$

In a Euclidean triangle with inner angles
$(\varphi(f,e),\varphi(f',e),\pi-D(e)),$ if the edge length $r(f)$ decreases and
the edge length $r(f')$ increases or does not change, the angle $\varphi(f,e)$
increases. Thus if $e\in E(F_0)_1$ and $f\in F_0$, as a function of $r$, $\varphi(f,e)$ satisfies
$\varphi(f,e)(r)<\varphi(f,e)(\overline{r}).$ Therefore $\sum_{f\in
F_0}\frac12\Phi(f)(r)<\sum_{f\in F_0}\frac12\Phi(f)(\overline{r}).$ This
shows $P:\Delta\to\delta$ is one to one.

Therefore $P:\Delta\to\delta$ is a smooth embedding. Thus
$P(\Delta)$ is open in $\delta.$ We will show $P(\Delta)$ is also
closed in $\delta.$ This will conclude that $P(\Delta)=\delta$
since $\delta$ is connected. To show the closeness we only need to
show that $P$ sends the boundary of $\overline{\Delta}$ to the
boundary of $\overline{\delta}.$ More precisely, if there is a
sequence of points $r^{(m)}$ in $\Delta$ so that
$\lim_{m\to\infty}r^{(m)}=s=(s(f_1),...,s(f_n))\in\partial\overline{\Delta},$
then $\lim_{m\to\infty}P(r^{(m)})\in\partial\overline{\delta},$
where we denote the boundary of a closed set $A$ by $\partial A.$
Since $s\in\partial\overline{\Delta},$ $F_1=\{f_i\in F|s(f_i)=0\}$
is a nonempty proper subset. We have
$$\sum_{f\in F_1}\frac12\Phi(f)(r^{(m)})=\sum_{e\in E(F_1)_1,f\in F_1}\varphi(f,e)(r^{(m)})+\sum_{e\in E(F_1)_2} D(e).$$
If $e\in E(F_1)_1,e<f\in F_1$ and $e<f'\notin F_1$, since $\lim_{m\to\infty}(r(f)(r^{(m)}),r(f')(r^{(m)}))=(s(f)=0,s(f')>0),$
we see that
$\lim_{m\to\infty}\varphi(f,e)(r^{(m)})=D(e).$ Hence
\begin{equation*}
\begin{split}
\lim_{m\to\infty}\sum_{f\in F_1}\frac12\Phi(f)(r^{(m)})
&=\sum_{e\in E(F_1)_1}D(e)+\sum_{e\in E(F_1)_2} D(e)\\
&=\sum_{e\in E(F_1)} D(e).
\end{split}
\end{equation*}
By the definition of $\delta,$ this shows that
$\lim_{m\to\infty}P(r^{(m)})\in\partial\overline{\delta}$.
\end{proof}

\section{Continuity method for hyperbolic circle pattern}

\begin{lemma}\label{degenerate} For a hyperbolic triangle with edge lengths $l,r_1,r_2$ and opposite angles
$\theta, \varphi_1, \varphi_2$ respectively, when $\theta\in
(0,\pi)$ is fixed, $ \varphi_1, \varphi_2$ are functions of $r_1,
r_2$.
\begin{enumerate}
\item If $r_1$ decreases and $r_2$ increases or does not change,
then $\varphi_2$ increases.

\item If both of $r_1, r_2$ decrease, then $\varphi_1 + \varphi_2$
increase.

\item If $(r_1,r_2)$ converges to $(0,a)$, where $a\in(0,\infty)$,
then $\varphi_2$ converges to $\pi-\theta.$

\item If $(r_1,r_2)$ converges to $(0,0)$, then $\varphi_1 +
\varphi_2$ converges to $\pi-\theta.$

\item If $(r_1,r_2)$ converges to $(\infty,b)$, where
$b\in[0,\infty)$, then $\varphi_2$ converges to $0.$

\item If $(r_1,r_2)$ converges to $(\infty,\infty),$ then
$\varphi_2$ converges to to $0.$
\end{enumerate}
\end{lemma}

\begin{figure}[htbp]
\begin{center}
\includegraphics[scale=.65]{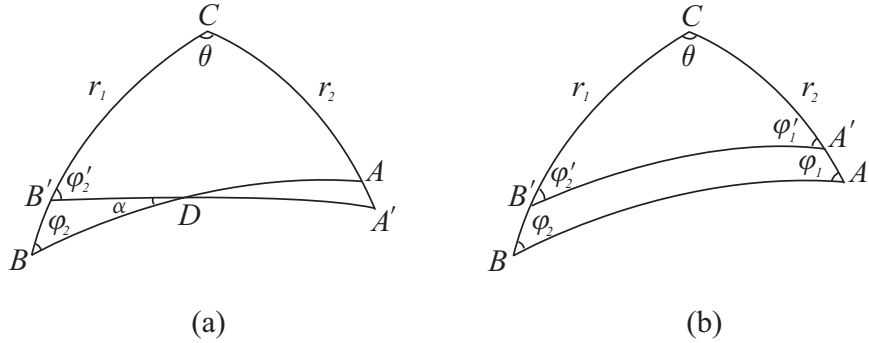}
\end{center}
\caption{\label{triangle}Angles change according to edge lengths}
\end{figure}

\begin{proof} (1) As in the Figure \ref{triangle} (a), $\triangle_{ABC}$ is the
hyperbolic triangle with edge lengths $r_1,r_2$ and angles $\varphi_1, \varphi_2.$
If $r_1$ decreases and $r_2$ increases or does not change, we obtain another hyperbolic triangle $\triangle_{A'B'C}$
with angles $\varphi_1', \varphi_2'.$ Consider the hyperbolic triangle $\triangle_{BDB'}$, we have
$\alpha+\varphi_2+\pi-\varphi_2'<\pi.$ Thus
$\varphi_2<\alpha+\varphi_2<\varphi_2'.$ This shows that $\varphi_2$ increases.

(2) As in the Figure \ref{triangle} (b), if both of $r_1, r_2$
decrease, we obtain another hyperbolic triangle $\triangle_{A'B'C}$ with angles $\varphi_1', \varphi_2'.$
Since the area of hyperbolic triangle $\triangle_{A'B'C}$ is
less than that of $\triangle_{ABC}$, we see
$\pi-(\theta+\varphi_1' + \varphi_2')<\pi-(\theta+\varphi_1 +
\varphi_2).$ Thus $\varphi_1 + \varphi_2<\varphi_1' + \varphi_2'.$

(3) If $\lim(r_1,r_2)=(0,a)$ with $a\in(0,\infty),$ by the cosine
law, we see
\begin{equation*}
\begin{split}
\lim\cosh l
&=\lim(\cosh r_1\cosh r_2-\cos \theta\sinh r_1 \sinh r_2)\\
&=\cosh a.
\end{split}
\end{equation*}
Thus $\lim l=a.$ By the cosine law, we see
\begin{equation*}
\begin{split}
\lim\cos\varphi_1&=\lim\frac{-\cosh r_1+\cosh r_2 \cosh l}{\sinh
r_2 \sinh l}\\
&=\frac{-1+\cosh^2 a}{\sinh^2 a}\\
&=1.
\end{split}
\end{equation*}
Thus $\lim\varphi_1=0.$ Since the area of the hyperbolic triangle
converges to $0.$ We have $\lim\varphi_2=\pi-\theta.$

(4) If $\lim(r_1,r_2)=(0,0),$ the area of the hyperbolic
triangle converges to $0,$ then $\lim(\varphi_1+\varphi_2)=
\pi-\theta.$

(5) If $\lim(r_1,r_2)=(\infty,b)$ with $b\in[0,\infty),$ by the
triangle inequality $l>r_1-r_2,$ then $\lim l=\infty.$ Thus
\begin{equation*}
\begin{split}
\lim\cos\varphi_2
&=\lim\frac{-\cosh r_2+\cosh r_1 \cosh l}{\sinh r_1 \sinh l}\\
&=\lim\frac{-\cosh r_2}{\sinh r_1 \sinh l}+\lim\frac{\cosh r_1}{\sinh r_1}\frac{\cosh l}{\sinh l}\\
&=0+1.
\end{split}
\end{equation*}
Hence $\lim\varphi_2=0.$

(6) If $\lim(r_1,r_2)=(\infty,\infty),$ by the cosine law, we see $$\lim\cosh l=\lim(\cosh r_1\cosh r_2-\cos
\theta\sinh r_1 \sinh r_2)=\lim e^{r_1+r_2}(1-\cos\theta)=\infty.$$ Thus $\lim l=\infty.$ By the cosine law,
\begin{equation*}
\begin{split}
\lim\cos\varphi_2
&=\lim\frac{-\cosh r_2+\cosh r_1 \cosh l}{\sinh r_1 \sinh l}\\
&=\lim\frac{-\cosh r_2}{\sinh r_1 \sinh l}+1\\
&=\lim\frac{-\cosh r_2}{\sinh r_1 \cosh l}+1\\
&=\lim\frac{-e^{r_2}}{e^{r_1}e^{r_1+r_2}(1-\cos\theta)}+1\\
&=0+1.
\end{split}
\end{equation*}
Hence $\lim\varphi_2=0.$
\end{proof}

\begin{proof}[Proof of Theorem \ref{h}] Assume the set of faces of the cell
decomposition is $F=\{f_1,...,f_n\}$. Consider the space of radii
$\Delta_H=\{(r(f_1),...,r(f_n))\}=\mathbb{R}^n_{>0}.$ For a fixed
intersection angle function $D:E\to(0,\pi)$, as in the Euclidean
case there is a mapping $P_H:\Delta_H\to \mathbb{R}^n$ sending the
radii $(r(f_1),...,r(f_n))$ to its cone angles
$\Phi=(\Phi(f_1),...,\Phi(f_n)).$ Now we introduce the polytope
$\delta_H=\{\Phi\in \mathbb{R}^n_{>0}|\Phi\ \text{satisfies the
condition}\ (\ref{hh})\}.$ We claim that $P_H(\Delta_H)\subset
\delta_H.$

In fact, for any nonempty subset $X$ of $F$, let's recall the decomposition $E(X)=E(X)_1\cup E(X)_2.$
Notice as in Figure \ref{Fig:circle}, in each hyperbolic triangle we have $\varphi(f,e)+\varphi(f',e)<D(e).$
Thus
\begin{equation*}
\begin{split}
\sum_{f\in X}\frac12\Phi(f) &=\sum_{f\in
X}\sum_{e<f}\varphi(f,e)\\
&=\sum_{e\in E(X)_1,f\in X}\varphi(f,e)+\sum_{e\in
E(X)_2}(\varphi(f,e)+\varphi(f',e))\\
&<\sum_{e\in E(X)_1,f\in X}\varphi(f,e)+\sum_{e\in E(X)_2} D(e)\\
&<\sum_{e\in E(X)_1}D(e)+\sum_{e\in E(X)_2} D(e)\\
&=\sum_{e\in E(X)}D(e).
\end{split}
\end{equation*}
Next, we claim the mapping $P_H:\Delta_H\to\delta_H$ is one to
one. In fact, let $r=(r(f_1),...,r(f_n)),
\overline{r}=(\overline{r}(f_1),...,\overline{r}(f_n))$ be two
distinct points in $\Delta_H.$ Let $F_0=\{f_i\in
F|r(f_i)>\overline{r}(f_i)\}$ and $\overline{F}_0=\{f_i\in
F|r(f_i)<\overline{r}(f_i)\}$. Since $r\neq r',$ at least one of
$F_0, \overline{F}_0$ is nonempty, say $F_0$ is nonempty. Now we
have
$$\sum_{f\in F_0}\frac12\Phi(f)=\sum_{e\in E(F_0)_1,f\in F_0}\varphi(f,e)+\sum_{e\in E(F_0)_2}(\varphi(f,e)+\varphi(f',e)).$$

If $e\in E(F_0)_1,e<f\in F_0$ and $e<f'\notin F_0$, then $r(f)$ decreases while $r(f')$ increases or does not change,
by Lemma \ref{degenerate}(1),
$\varphi(f,e)(r)<\varphi(f,e)(\overline{r}).$ If $e\in E(F_0)_2$, then both
of $r(f)$ and $r(f')$ decrease, by Lemma \ref{degenerate}(2),
$\varphi(f,e)(r)+\varphi(f',e)(r)<\varphi(f,e)(\overline{r})+\varphi(f',e)(\overline{r}).$
Therefore $\sum_{f\in F_0}\frac12\Phi(f)(r)<\sum_{f\in
F_0}\frac12\Phi(f)(\overline{r}).$ This shows $P_H:\Delta_H\to\delta_H$ is
one to one.

Therefore Therefore $P_H:\Delta_H\to\delta_H$ is a smooth
embedding. Thus $P_H(\Delta_H)$ is open in $\delta_H.$ We will
show $P_H(\Delta_H)$ is also closed in $\delta_H.$ This will
conclude that $P_H(\Delta_H)= \delta_H$ since $\delta_H$ is
connected. To show the closeness, we assume $r^{(m)}$ is a
sequence of points in $\Delta_H$ such that $\lim_{m\to\infty}
P_H(r^{(m)})=t\in\delta_H.$ By taking subsequence, we may assume
$\lim_{m\to\infty} r^{(m)}=s\in[0,\infty]^n.$ We only need to show
$s\in(0,\infty)^n=\Delta_H.$ This will finish the proof since
$P_H(s)=t.$

Suppose otherwise that there is some face $f$ so that
$s(f)\in\{0,\infty\},$ we will discuss the two cases.

Case 1, if there is some face $f$ so that $\lim_{m\to\infty}r^{(m)}(f)=s(f)=\infty,$ by Lemma
\ref{degenerate}(5) and (6), we have
$\lim_{m\to\infty}\varphi(f,e)(r^{(m)})=0$ for each $e<f.$ Thus
\begin{equation*}
\begin{split}
t(f)&=\lim_{m\to\infty}\Phi(f)(r^{(m)})\\
&=\lim_{m\to\infty}\sum_{e<f}\varphi(f,e)(r^{(m)})\\
&=0.
\end{split}
\end{equation*}
This is a contradiction since we assume $t\in\delta_H.$

Case 2, there is no face $f$ such that $s(f)=\infty$ but the
subset $F_1=\{f_i\in F|s(f_i)=0\}$ is nonempty. We have
$$\sum_{f\in F_1}\frac12\Phi(f)(r^{(m)})=\sum_{e\in E(F_1)_1,f\in F_1}\varphi(f,e)(r^{(m)})+\sum_{e\in E(F_1)_2} (\varphi(f,e)+\varphi(f',e))(r^{(m)}).$$
If $e\in E(F_1)_1,e<f\in F_1$ and $e<f'\notin F_1$, then
$\lim_{m\to\infty}(r^{(m)}(f),r^{(m)}(f'))=(0,s(f')>0)$, by Lemma
\ref{degenerate}(3),
$\lim_{m\to\infty}\varphi(f,e)(r^{(m)})=D(e).$ If $e\in E(F_1)_2$,
then $\lim_{m\to\infty}(r^{(m)}(f),r^{(m)}(f'))=(0,0),$ by Lemma
\ref{degenerate} (4),
$\lim_{m\to\infty}(\varphi(f,e)+\varphi(f',e))(r^{(m)})=D(e).$
Hence
\begin{equation*}
\begin{split}
\lim_{m\to \infty}\sum_{f\in
F_1}\frac12\Phi(f)(r^{(m)})
&=\sum_{e\in E(F_1)_1}D(e)+\sum_{e\in E(F_1)_2} D(e)\\
&=\sum_{e\in E(F_1)} D(e).
\end{split}
\end{equation*}
This is a contradiction since we assume $t\in\delta_H.$
\end{proof}

\section {Variational principle for hyperbolic circle pattern}

\subsection{Angle structure}
An \it angle structure \rm for a hyperbolic circle pattern with an
intersection angle function $D:E\to(0,\pi)$ and a cone angle
function $\Phi:F\to(0,\infty)$ is a function $\varphi: F\times
E\to (0,\pi),$ such that $\varphi(f,e)+\varphi(f',e)<D(e)$ for
each two faces $f,f'$ sharing an edge $e$ and
$\sum_{e<f}\varphi(f,e)=\frac12\Phi(f)$ for each face $f$. Denote
the space of all angle structures by $\mathcal{A}.$ In this
section we will use the duality theorem in the linear programming.
For a statement and a proof of the theorem, see for example Kolman
and Beck \cite{kb}.

\begin{lemma} \label{l} Given a cellular decomposed surface
$\Sigma$, an intersection angle function $D:E\to (0,\pi)$ and a
cone angle function $\Phi: F \to (0,\infty)$, the closure of
$\mathcal{A}$ is nonempty if and only if
\begin{equation}
\label{l} \sum_{f\in X}\frac12\Phi(f)\leq \sum _{e\in E(X)} D(e)
\end{equation}

for any $X \subseteq F.$
\end{lemma}

\begin{proof}Consider the linear programming problem $(P)$: minimize the
objective function $z = 0$, subject to the constraints
$$
\left\{
\begin{array}{lll}
\varphi(f,e)+\varphi(f',e)\leq D(e)\\
\sum_{e<f}\varphi(f,e)=\frac12\Phi(f)\\
\varphi(f,e)\geq 0
\end{array}
\right.
$$

The closure of $\mathcal{A}$ is exactly the set of solutions of
above inequalities.

The dual problem $(P^*)$ has variables $( ...,y_f, ..., y_e, ...)$
indexed by $F \cup E$. The dual problem $(P^*)$ is to maximize the
objective function $z = \sum_{f \in F} \frac12\Phi(f) y_f+ \sum_{e
\in E} D(e) y_e,$ subject to the constraints
$$
\left\{
\begin{array}{ll}
y_e\leq0\\
y_f + y_e \leq 0\ \ \mbox{whenever}\ e < f
\end{array}
\right.
$$

Now by the definition of $\mathcal{A}$ and the duality theorem in the linear programming,
the following statements are equivalent.
\begin{enumerate}
\item[(i)] The closure of $\mathcal{A}$ is nonempty. \item[(ii)]
There exists a feasible solution of the problem $(P)$.
\item[(iii)] The objective function of the problem $(P^*)$ is
non-positive.
\end{enumerate}
Hence to prove this lemma we only need to show the condition
(\ref{l}) is equivalent to the statement (iii).

To show the condition (\ref{l}) is necessary, we assume that the objective function of $(P^*)$ is
non-positive. For any $X \subseteq
F,$ let
$$y_f=
\left\{
\begin{array}{ccc}
1&\mbox{if}\ f\in X \\
0&\mbox{if}\ f\notin X
\end{array}
\right.\     \mbox{and}\     y_e= \left\{
\begin{array}{ccc}
-1&\mbox{if}\ e\in E(X) \\
0&\mbox{if}\ e\notin E(X)
\end{array}
\right.
$$
It is easy to check $(y_f,y_e)$ is a feasible solution of $(P^*)$. Thus
$0\geq z(y_f, y_e) = \sum_{f \in X} \frac12\Phi(f) y_f+ \sum_{e
\in E(X)} D(e) y_e = \sum_{f\in X}\frac12\Phi(f)-\sum _{e\in
E(X)}D(e).$ This gives the condition (\ref{l}).

To show the condition (\ref{l}) is sufficient, take an arbitrary
feasible solution of $(P^*)$: $(y_f, y_e)$. Since we can replace
the negative $y_f$'s by $0$'s such that the objective function
does not decrease, we can assume all $y_f\geq 0.$ Let $X=\{f\in F
| y_f > 0\}$, and $a$ be the minimal element in the set $\{y_f,
f\in X\}.$ Thus $a>0.$ Then we define
$$y_f^{(1)}=
\left\{
\begin{array}{ccc}
y_f -a&\mbox{if}\ f\in X \\
y_f =0&\mbox{if}\ f\notin X
\end{array}
\right.\     \mbox{and}\     y_e^{(1)}= \left\{
\begin{array}{ccc}
y_e +a&\mbox{if}\ e\in E(X) \\
y_e &\mbox{if}\ e\notin E(X)
\end{array}
\right.
$$
We can check $(y_f^{(1)},y_e^{(1)})$ is still a feasible solution of $(P^*)$. In
fact if $e\notin E(X),$ then $y_e^{(1)}=y_e\leq0.$ If $e\in E(X),$
there exists $f'\in X$ such that $e<f'.$ Then $y_e^{(1)}=y_e
+a\leq y_e +y_{f'}\leq 0,$ since $a$ is minimal and by assumption $(y_{f'},y_e)$ is
feasible. Next we consider $y_f^{(1)}+y_e^{(1)}.$ If $f\in X,$
$y_f^{(1)}+y_e^{(1)}=y_f+y_e\leq 0$ or $y_f^{(1)}+y_e^{(1)}=y_f-a+y_e<0.$ If $f\notin
X,$ $y_f^{(1)}+y_e^{(1)}=y_e^{(1)}\leq 0.$ Hence
$(y_f^{(1)},y_e^{(1)})$ is feasible.

Now $z(y_f^{(1)},y_e^{(1)})=z(y_f, y_e)+a(\sum _{e\in
E(X)}D(e)-\sum_{f\in X}\frac12\Phi(f))\geq z(y_f, y_e).$ due to
the condition (\ref{l}). Notice that the number of 0's in
$\{y_f^{(1)}\}$ is more than that in $\{y_f \}$. By the same
procedure, we obtain a sequence of feasible solutions
$(y_f^{(i)},y_e^{(i)})$. After finite steps, it ends at a feasible
solution $(y_f^{(n)}=0,y_e^{(n)}\leq0)$. Thus $z(y_f, y_e)\leq
z(y_f^{(1)},y_e^{(1)})\leq ...\leq
z(y_f^{(n)},y_e^{(n)})=\sum_{e\in E}D(e)y_e^{(n)}\leq0.$ This
shows that the objective function of $(P^*)$ is non-positive.
\end{proof}

\begin{theorem}[Bobenko-Springborn \cite{bs} Proposition 4]\label{angle} Given a cellular decomposed surface
$\Sigma$, an intersection angle function $D:E\to (0,\pi)$ and a
cone angle function $\Phi: F \to (0,\infty)$, $\mathcal{A}$ is
nonempty if and only if

\begin{equation}
\label{ll} \sum_{f\in X}\frac12\Phi(f)< \sum _{e\in E(X)} D(e)
\end{equation}

for any $X \subseteq F.$
\end{theorem}

\begin{proof} Let
$\varphi(f,e)=a(f,e)+\varepsilon,$ where $a(f,e)\geq 0$ and
$\varepsilon\geq 0.$ Consider the linear programming problem $(P)$
with variables $\{...,a(f,e),...\varepsilon\}$: minimize the
objective function $z = -\varepsilon$, subject to the constraints
$$
\left\{
\begin{array}{llll}
a(f,e)+a(f',e)+2\varepsilon \leq D(e)\\
\sum_{e<f}a(f,e)+n_f\varepsilon=\frac12\Phi(f)\\
a(f,e)\geq 0 \\
\varepsilon\geq 0
\end{array}
\right.
$$
where $n_f$ is the number of edges of the face $f.$

The dual problem $(P^*)$ has variables $( ...,y_f, ..., y_e, ...)$
indexed by $F \cup E$. The dual problem $(P^*)$ is to maximize the
objective function $z = \sum_{f \in F} \frac12\Phi(f) y_f+ \sum_{e
\in E} D(e) y_e,$ subject to the constraints
$$
\left\{
\begin{array}{lll}
y_e\leq0\\
y_f + y_e \leq 0\ \  \  \  \ \ \  \mbox{whenever}\ e < f \\
\sum_{f\in F}n_fy_f+2\sum_{e\in E}y_e \leq -1
\end{array}
\right.
$$
Now by the definition of $\mathcal{A}$ and the duality theorem in the linear programming,
the following statements are equivalent.
\begin{enumerate}
\item[(i)] The set $\mathcal{A}$ is nonempty. \item[(ii)] There
exists a feasible solution of the problem $(P)$ such that
$a(f,e)\geq 0, \varepsilon> 0$. \item[(iii)] The minimal value of
the objective function of the problem $(P)$ is negative.
\item[(iv)] The objective function of the problem $(P^*)$ is
negative.
\end{enumerate}
Hence to prove this lemma we only need to show the condition
(\ref{ll}) is equivalent to the statement (iv).

To show the condition (\ref{ll}) is necessary, if $X=F,$ then
\begin{equation*}
\begin{split}
\sum_{f\in F}\frac12\Phi(f) &=\sum_{f\in
F}\sum_{e<f}\varphi(f,e)\\
&=\sum_{e\in E}(\varphi(f,e)+\varphi(f',e))\\
&<\sum_{e\in E}D(e).
\end{split}
\end{equation*}
For any $X \subset F,X\neq F$ let
$$y_f=
\left\{
\begin{array}{ccc}
1&\mbox{if}\ f\in X \\
0&\mbox{if}\ f\notin X
\end{array}
\right.\     \mbox{and}\     y_e= \left\{
\begin{array}{ccc}
-1&\mbox{if}\ e\in E(X) \\
0&\mbox{if}\ e\notin E(X)
\end{array}
\right.
$$
To show $(y_f,y_e)$ is a feasible solution of $(P^*)$ we only need
to check it satisfies $\sum_{f\in F}n_fy_f+2\sum_{e\in E}y_e \leq
-1.$ This is equivalent to $\sum_{f\in X}n_f -2\sum_{e\in E(X)}
\leq -1.$ It is true since $X\neq F$ implies that
$2|E(X)|>\sum_{f\in X}n_f$ or $2|E(X)|\geq\sum_{f\in X}n_f+1.$
Since the objective function is negative, we have $0>z(y_f,y_e)$
which implies $\sum_{f\in X}\frac12\Phi(f)< \sum _{e\in E} D(e).$
This gives the condition (\ref{ll}).

To show the condition (\ref{ll}) is sufficient, by the proof of Lemma
\ref{l} we know the objective function of
$(P^*)$ is $\leq 0$ under the condition (\ref{ll}) which is stronger than the condition (\ref{l}). We try to show it can
not be 0. Otherwise, let $(y_f,y_e)$ be a feasible solution satisfying
$z(y_f,y_e)=0.$ We obtain a new feasible solution $(y'_f,y_e)$ of $(P^*)$ by replacing
the negative $y_f$'s by $0$'s. Since the objective function does
not decrease, $0=z(y_f,y_e)\leq z(y'_f,y_e)\leq 0$. If there is some face $f$ such that $y'_f>0,$ we can construct
another feasible solution $(y_f^{(1)},y_e^{(1)})$ as in the proof of Lemma \ref{l} so that
$z(y_f^{(1)},y_e^{(1)})>z(y'_f,y_e)=0.$ This is impossible since the objective function $z\leq 0$.

Hence $y'_f=0$ for each $f.$ Therefore $0=z(y'_f,y_e)=\sum _{e\in
E} D(e)y_e.$ Since $y_e\leq 0$ for each $e,$ we see that $y_e=0$
for each $e.$ But $(y'_f=0,y_e=0)$ does not satisfy the condition
$\sum_{f\in F}n_fy_f+2\sum_{e\in E}y_e \leq -1$. This
contradiction implies that the objective function of $(P^*)$ can
not be 0.
\end{proof}

\subsection{The derivative cosine law}
To prove Theorem \ref{h} we will define an energy function on
$\mathcal{A}$ the space of angle structures. First, we will
construct a function for one triangle.

The following lemma is essentially obtained in Luo \cite{l1}. For the convenience of
readers, we include a simple proof here. We use
$$\left[
\begin{array}{ccc}
a_1\\
\vdots\\
a_n
\end{array}
\right]$$ to denote the diagonal matrix with diagonal entries
$(a_1,\ldots,a_n)$.
\begin{lemma}[The derivative cosine law]\label{derivative} For a hyperbolic triangle with three edges
$l_1,l_2,l_3$ and opposite angles $\theta_1,\theta_2,\theta_3,$
$\{i,j,k\}=\{1,2,3\}$, the differentials of $l'$s and $\theta'$s satisfy
$$\left(
\begin{array}{ccc}
dl_1 \\
dl_2 \\
dl_3
\end{array}\right)
=\frac{-1}{\sinh l_i\sinh l_j\sin \theta_k} \left[
\begin{array}{ccc}
\sinh l_1\\
\sinh l_2\\
\sinh l_3
\end{array}
\right] \left(
\begin{array}{ccc}
1&\cosh l_3&\cosh l_2 \\
\cosh l_3&1&\cosh l_1 \\
\cosh l_2&\cosh l_1&1
\end{array}
\right) \left(
\begin{array}{ccc}
d\theta_1 \\
d\theta_2 \\
d\theta_3
\end{array}\right).$$
\end{lemma}
\begin{proof} By differentiating the cosine law of hyperbolic triangles
$$\cosh l_i\sin\theta_j\sin \theta_k=\cos \theta_i+\cos \theta_j\cos \theta_k,$$
we have
$$\sinh l_i\sin \theta_j\sin \theta_kdl_i$$$$=-\sin \theta_id\theta_i
-(\sin \theta_j\cos \theta_k+\cosh l_i\cos\theta_j\sin
\theta_k)d\theta_j -(\cos \theta_j\sin \theta_k+\cosh l_i\sin
\theta_j\cos \theta_k)d\theta_k.$$ After replacing $\cosh l_i$ by
using the cosine law $$\cosh l_i=\frac{\cos \theta_i+\cos
\theta_j\cos \theta_k}{\sin\theta_j\sin \theta_k}$$ and
simplifying we get
\begin{equation*}
\begin{split}
\sinh l_i\sin\theta_j\sin \theta_kdl_i &=-\sin \theta_id\theta_i-
\frac{\cos \theta_i\cos\theta_j+\cos \theta_k}{\sin
\theta_j}d\theta_j-\frac{\cos \theta_i\cos\theta_k+\cos
\theta_j}{\sin \theta_k}d\theta_k\\
&=-\sin \theta_i(d\theta_i+\cosh l_kd\theta_j+\cosh l_jd\theta_k).
\end{split}
\end{equation*}
The matrix is derived from the sine law:
$$\frac{\sin\theta_i}{\sinh l_i}=\frac{\sin\theta_j}{\sinh l_j},\ \text{for any}\ i,j\in\{1,2,3\}.$$
\end{proof}

Let the domain
$T=\{(\varphi_1,\varphi_2)|\varphi_1>0,\varphi_2>0,\varphi_1+\varphi_2<\pi-\theta\}$
be the space of angles $(\varphi_1,\varphi_2)$ in a hyperbolic
triangle with the third angle $\theta$ fixed.
\begin{lemma}\label{cosine}
For a hyperbolic triangle with three edges $r_1,r_2,l$ and
opposite angles $\varphi_1, \varphi_2,\theta,$ when $\theta\in
(0,\pi)$ is fixed, let $w=\ln\tanh\frac{r_1}2 d\varphi_2
+\ln\tanh\frac{r_2}2 d\varphi_1 $. Then the function
$v(\varphi_1,\varphi_2)=\int_{(0,0)}^{(\varphi_1,\varphi_2)} w$ is
well defined and strictly concave down in the domain $T.$
\end{lemma}
\begin{proof} If $\theta\in
(0,\pi)$ is fixed, then $d\theta=0,$ by Lemma \ref{derivative}, we have
\begin{equation*}
\begin{split}
\left(
\begin{array}{ccc}
dr_1 \\
dr_2
\end{array}\right)
&=\frac{-1}{\sinh r_1\sinh r_2\sin \theta} \left[
\begin{array}{ccc}
\sinh r_1\\
\sinh r_2
\end{array}
\right] \left(
\begin{array}{ccc}
1&\cosh l \\
\cosh l&1
\end{array}
\right) \left(
\begin{array}{ccc}
d\varphi_1 \\
d\varphi_2
\end{array}\right)\\
&=\frac{-1}{\sinh r_1\sinh r_2\sin \theta} \left[
\begin{array}{ccc}
\sinh r_1\\
\sinh r_2
\end{array}
\right] \left(
\begin{array}{ccc}
\cosh l & 1 \\
1 & \cosh l
\end{array}
\right) \left(
\begin{array}{ccc}
d\varphi_2 \\
d\varphi_1
\end{array}\right).
\end{split}
\end{equation*}
Then
\begin{equation*}
\begin{split}
\left(
\begin{array}{ccc}
d\ln\tanh\frac{r_1}2 \\
d\ln\tanh\frac{r_2}2
\end{array}\right)
&=\left[
\begin{array}{ccc}
\sinh^{-1} r_1\\
\sinh^{-1} r_2
\end{array}
\right]\left(
\begin{array}{ccc}
dr_1 \\
dr_2
\end{array}\right)\\
&=\frac{-1}{\sinh r_1\sinh r_2\sin \theta} \left(
\begin{array}{ccc}
\cosh l & 1 \\
1 & \cosh l
\end{array}
\right) \left(
\begin{array}{ccc}
d\varphi_2 \\
d\varphi_1
\end{array}\right).
\end{split}
\end{equation*}
We denote it to be
$$\left(
\begin{array}{ccc}
d\ln\tanh\frac{r_1}2 \\
d\ln\tanh\frac{r_2}2
\end{array}\right)
=M\left(
\begin{array}{ccc}
d\varphi_2 \\
d\varphi_1
\end{array}\right),$$
where $M$ is symmetric and negative definite. Hence
$w=\ln\tanh\frac{r_1}2 d\varphi_2 +\ln\tanh\frac{r_2}2 d\varphi_1
$ is a closed $1-$form in the domain $T$. Therefore the function
$v(\varphi_1,\varphi_2)$ is well defined. Since the Hessian matrix
of $v(\varphi_1,\varphi_2)$ is exactly the matrix $M$ which is
negative definite, $v(\varphi_1,\varphi_2)$ is strictly concave
down in $T.$
\end{proof}

\subsection{Boundary behavior}
In this section we will investigate the behavior of the function
$v(\varphi_1,\varphi_2)$ on the boundary of the closure of the
domain $T.$
\begin{lemma}\label{extend} The function $v(\varphi_1,\varphi_2)$ can be continuously extended to
the boundary of the closure of the domain $T.$
\end{lemma}
\begin{proof}
By the cosine law
\begin{equation*}
\begin{split}
(\tanh\frac{r_1}2)^2&=\frac{\cosh r_1-1}{\cosh r_1+1}\\
&=\frac{\cos\varphi_1+\cos(\varphi_2+\theta)}{\cos\varphi_1-\cos(\varphi_2+\theta)}\\
&=\frac{\sin\frac{\pi-\theta-\varphi_1-\varphi_2}2\sin\frac{\pi+\theta-\varphi_1+\varphi_2}2}
{\sin\frac{\pi+\theta-\varphi_1-\varphi_2}2\sin\frac{\pi-\theta-\varphi_1+\varphi_2}2}.
\end{split}
\end{equation*}
By interchange the index 1,2, we obtain the similar formula.
\begin{equation*}
\begin{split}
(\tanh\frac{r_2}2)^2&=\frac{\sin\frac{\pi-\theta-\varphi_1-\varphi_2}2\sin\frac{\pi+\theta-\varphi_2+\varphi_1}2}
{\sin\frac{\pi+\theta-\varphi_1-\varphi_2}2\sin\frac{\pi-\theta-\varphi_2+\varphi_1}2}\\
&=\frac{\sin\frac{\pi-\theta-\varphi_1-\varphi_2}2\sin\frac{\pi-\theta-\varphi_1+\varphi_2}2}
{\sin\frac{\pi+\theta-\varphi_1-\varphi_2}2\sin\frac{\pi+\theta-\varphi_1+\varphi_2}2}
\end{split}
\end{equation*}
due to
$\frac{\pi+\theta-\varphi_2+\varphi_1}2+\frac{\pi-\theta-\varphi_1+\varphi_2}2=\pi.$

Thus
\begin{equation*}
\begin{split}
w&=\frac12(\ln\sin\frac{\pi-\theta-\varphi_1-\varphi_2}2-\ln\sin\frac{\pi+\theta-\varphi_1-\varphi_2}2)(d\varphi_2+d\varphi_1)\\
&+\frac12(\ln\sin\frac{\pi+\theta-\varphi_1+\varphi_2}2-
\ln\sin\frac{\pi-\theta-\varphi_1+\varphi_2}2)(d\varphi_2-d\varphi_1).
\end{split}
\end{equation*}

Hence
\begin{equation*}
\begin{split}
v(\varphi_1,\varphi_2)&=\Lambda(\frac{\pi-\theta-\varphi_1-\varphi_2}2)+
\Lambda(\frac{\pi-\theta+\varphi_1+\varphi_2}2)\\
&+\Lambda(\frac{\pi-\theta+\varphi_1-\varphi_2}2)+
\Lambda(\frac{\pi-\theta-\varphi_1+\varphi_2}2)-4\Lambda(\frac{\pi-\theta}2),
\end{split}
\end{equation*}
where $\Lambda(x)=- \int_0^xln|2\sin t|dt$ is the Lobachevsky
function. By the property of Lobachevsky function, function
$v(\varphi_1,\varphi_2)$ can be continuously extended to the
boundary of the domain $T$.
\end{proof}

The closure of the domain
$T=\{(\varphi_1,\varphi_2)|\varphi_1>0,\varphi_2>0,\varphi_1+\varphi_2<\pi-\theta\}$
is a triangle with vertexes $O=(0,0), V_1=(\pi-\theta,0),
V_2=(0,\pi-\theta).$ If $(\varphi_1,\varphi_2)$ is in the line
segment of $V_1V_2$, the triangle with angles
$(\varphi_1,\varphi_2,\theta)$ is a Euclidean or degenerated
Euclidean triangle. Otherwise the triangle with angles
$(\varphi_1,\varphi_2,\theta)$ is a hyperbolic or degenerated
hyperbolic triangle. Let $(\varphi_1(t),\varphi_2(t))$ be a line
segment, where, for $i=1,2$ and $0\leq t \leq 1,$
$$\varphi_i(t)=(1-t)\varphi_i(0)+t\varphi_i(1),$$
such that $(\varphi_1(1),\varphi_2(1))\in T$ and
$(\varphi_1(0),\varphi_2(0))$ is in the closure of $T$.

\begin{lemma} \label{b}
Let $v(t)=v(\varphi_1(t),\varphi_2(t))$.
\begin{enumerate}
\item If $(\varphi_1(0),\varphi_2(0))$ is in the line segment of
$V_1V_2$, then $\lim_{t\to 0}v'(t)=\infty$.

\item Otherwise, $\lim_{t\to 0}v'(t)$ exists and is a finite
number.
\end{enumerate}
\end{lemma}

\begin{proof} We have $$v'(t)=\frac
d{dt}v(\varphi_1(t),\varphi_2(t)) =\varphi_2'(t)\ln\tanh
\frac{r_1(t)}2+ \varphi_1'(t)\ln\tanh \frac{r_2(t)}2,$$ where $r_1(t)$ and $r_2(t)$
are the edge lengths of the hyperbolic triangle with angles $(\varphi_1(t),\varphi_2(t),\theta).$

Case 1, if $(\varphi_1(0),\varphi_2(0))\in T$, in the hyperbolic
triangle with angles $(\varphi_1(0),\varphi_2(0), \theta)$, the
edge lengths $r_1(0), r_2(0)$ are finite. Hence $\lim_{t\to
0}v'(t)$ exists and is a finite number.

Case 2, if $(\varphi_1(0),\varphi_2(0))$ is in the interior of the
segment $OV_1$ (similar for $OV_2$), in the hyperbolic triangle
with angles $(\varphi_1(0),\varphi_2(0)=0, \theta)$, the edge
lengths $r_1(0)$ is infinite and the edge length $r_2(0)$ is
finite. Hence $\lim_{t\to 0}v'(t)=\varphi_1'(0)\ln\tanh
\frac{r_2(0)}2=(\varphi_1(1)-\varphi_1(0))\ln\tanh
\frac{r_2(0)}2$ is a finite number.

Case 3, if $(\varphi_1(0),\varphi_2(0))=O,$ in the hyperbolic
triangle with angles $(\varphi_1(0)=0,\varphi_2(0)=0, \theta)$,
the two edge lengths $r_1(0), r_2(0)$ are infinite. Hence
$\lim_{t\to 0}v'(t)=0.$

If $(\varphi_1(0),\varphi_2(0))$ is in the segment of $V_1V_2$,
$(\varphi_1(0),\varphi_2(0), \theta)$ are the angles of a
Euclidean or degenerated Euclidean triangle. Thus the two edge
lengths $r_1(0), r_2(0)$ turn to zero.

By the computation in Lemma \ref{extend}, we have
\begin{equation*}
\begin{split}
v'(t)&=\frac12(\varphi_2'(t)+\varphi_1'(t))(\ln\sin\frac{\pi-\theta-\varphi_1-\varphi_2}2-\ln\sin\frac{\pi+\theta-\varphi_1-\varphi_2}2)\\
&+\frac12(\varphi_2'(t)-\varphi_1'(t))(\ln\sin\frac{\pi+\theta-\varphi_1+\varphi_2}2-
\ln\sin\frac{\pi-\theta-\varphi_1+\varphi_2}2).
\end{split}
\end{equation*}
Case 4, if $(\varphi_1(0),\varphi_2(0))$ is in the interior of the
segment $V_1V_2,$ then $\varphi_1(0)+\varphi_2(0)+\theta=\pi$ and
$\varphi_1(0)>0,\varphi_2(0)>0.$ Thus $\lim_{t\to
0}v'(t)=\frac12(\varphi_2'(0)+\varphi_1'(0))\ln\sin 0+$finite terms.

Since
$\varphi_1'(0)+\varphi_2'(0)=(\varphi_1(1)+\varphi_2(1))-(\varphi_1(0)+\varphi_2(0))<(\pi-\theta)-(\pi-\theta)=0,$
then $\lim_{t\to 0}v'(t)=\infty.$

Case 5, if $(\varphi_1(0),\varphi_2(0))=V_1$ (similar for $V_2$),
then $\varphi_1(0)=\pi-\theta,\varphi_2(0)=0.$ Thus $\lim_{t\to
0}v'(t)=\frac12(\varphi_2'(0)+\varphi_1'(0))(\ln\sin
0-\ln\sin\theta)+\frac12(\varphi_2'(0)-\varphi_1'(0))(\ln\sin\theta-\ln\sin
0)=\varphi_1'(0)(\ln\sin 0-\ln\sin\theta)=(\varphi_1(1)-\varphi_1(0))(\ln\sin 0-\ln\sin\theta)$. Since
$\varphi_1(0)=\pi-\theta>\varphi_1(1),$ then
$\lim_{t\to 0}v'(t)=\infty.$

\end{proof}

\subsection{Proof of Theorem \ref{h}}

Given a cellular decomposed surface $\Sigma$, an intersection
angle function $D:E\to(0,\pi)$ and a cone angle function
$\Phi:F\to(0,\infty),$ if there exists a hyperbolic circle
pattern, for each edge $e$ there is a quadrilateral which is a
union of two congruent hyperbolic triangles. The angles of these
triangles define an angle structure. By Theorem \ref{angle}, the
condition (\ref{hh}) holds.

On the other hand, if the condition (\ref{hh}) holds, by Theorem \ref{angle}, $\mathcal{A}$ the space
of all angle structures is nonempty.
Given an angle structure $\varphi\in\mathcal{A}$, for each edge $e$ shared by face $f$ and $f'$, there
is a topological quadrilateral. We realize the quadrilateral as a union of two congruent
hyperbolic triangles with angles $(\varphi(f,e),\varphi(f',e),\pi-D(e)).$ Let $r(f,e)$ and $r(f',e)$ be
the edge lengths of one of the two triangles. Since the two triangles are always congruent, we only need
to take care one of them.

For each hyperbolic triangle we have defined a function $v.$ Summation of the function $v$
over ``half" of all of triangles produces an \it energy function \rm $\mathcal{E}:\mathcal{A}\to \mathbb{R}.$
In fact for $\varphi\in\mathcal{A},$ let
\begin{equation*}
\begin{split}
\mathcal{E}(\varphi)&=\sum_{e\in E}v(\varphi(f,e),\varphi(f',e))\\
&=\sum_{e\in
E}\int_{(0,0)}^{(\varphi(f,e),\varphi(f',e))}(\ln\tanh
\frac{r(f,e)}2d\varphi(f,e)+\ln\tanh
\frac{r(f',e)}2d\varphi(f',e)),
\end{split}
\end{equation*}
where faces $f,f'$ share the edge $e,$ $r(f,e)$ and $r(f',e)$ are
the edge lengths of the hyperbolic triangle with angles
$(\varphi(f,e),\varphi(f',e),\pi-D(e)).$ By Lemma \ref{cosine} and
Lemma \ref{extend}, the function $\mathcal{E}$ is strictly concave
down in $\mathcal{A}$ and can be continuously extended to the
boundary of the closure of $\mathcal{A}$. Since the closure of
$\mathcal{A}$ is compact, $\mathcal{E}$ has the unique maximum
point $\varphi_0$ in closure of $\mathcal{A}$.

First we claim that under the angle structure $\varphi_0$, there is no Euclidean or degenerated Euclidean triangles.
Otherwise, choose any point $\varphi_1\in\mathcal{A}$ and join $\varphi_0$ and $\varphi_1$ by a line segment
$\varphi_t=(1-t)\varphi_0+t\varphi_1.$ Let $E=I\cup J,$ where each edge in $I$ corresponds to a hyperbolic or degenerate
hyperbolic triangle while each edge in $J$ corresponds to a Euclidean or degenerated Euclidean triangle.
Then
$$\mathcal{E}(\varphi_t)=\sum_{e\in I}v(\varphi_t(f,e),\varphi_t(f',e))+\sum_{e\in J}v(\varphi_t(f,e),\varphi_t(f',e)).$$
Therefore
$$\lim_{t\to 0}\frac d{dt}\mathcal{E}(\varphi_t)=
\sum_{e\in I}\lim_{t\to 0}\frac d{dt}v(\varphi_t(f,e),\varphi_t(f',e))+
\sum_{e\in J}\lim_{t\to 0}\frac d{dt}v(\varphi_t(f,e),\varphi_t(f',e)).$$
By Lemma \ref{b}, every term with index in $I$ is finite while every term with index in $J$ is $\infty.$
This is a contradiction. Since we assume $\varphi_0$ is the maximum point, $\lim_{t\to 0}\frac d{dt}\mathcal{E}(\varphi_t)$
should be non-positive.

Now we may assume under the angle structure $\varphi_0$, there is no Euclidean or degenerated Euclidean triangles.
Second we claim that under the angle structure $\varphi_0$, there is no zero angle.
Otherwise some face $f$ has a triangle $\triangle_e$ with angles
$(\varphi(f,e)=0, x, \pi-D(e)).$ Since the sum
$\sum_{e<f}\varphi(f,e)=\frac12\Phi(f)>0,$ face $f$ has another
triangle $\triangle_{e'}$ with angles $(\varphi(f,e')=y_1>0,y_2,\pi-D(e')).$
Now we keep all the angles in the triangles other than $\triangle_e$ and $\triangle_{e'}$ while
modify the angles of $\triangle_e$ to $(t, x, \pi-D(e))$ and modify the angles of $\triangle_{e'}$ to
$(y_1-t,y_2,\pi-D(e')).$ Since there is no Euclidean or degenerated Euclidean triangles in $\varphi_0$, then
$0+x+\pi-D(e)<\pi.$ We can make $t+x+\pi-D(e)<\pi$ for
sufficient small $t>0.$ And since $y_1>0,$ we can make $y_1-t>0$ for sufficient small $t>0.$
Therefore this angle assignment gives a point $\varphi_t$ in the closure of $\mathcal{A}.$
Now $$\lim_{t\to 0}\frac d{dt}\mathcal{E}(\varphi_t)=\lim_{t\to 0}\frac d{dt}v(t,x)+\lim_{t\to 0}\frac d{dt}v(y_1-t,y_2).$$
Since
$$\frac d{dt}v(\varphi_1(t),\varphi_2(t))
=\varphi_2'(t)\ln\tanh \frac{r_1(t)}2+ \varphi_1'(t)\ln\tanh
\frac{r_2(t)}2,$$ we have
$$\lim_{t\to 0}\frac d{dt}v(t,x)=\lim_{t\to 0} [x'\ln\tanh \frac{r_1(t)}2+ t'\ln\tanh \frac{r_2(t)}2]
=\lim_{t\to 0}\ln\tanh \frac{r_2(t)}2,$$ where $r_1(t)$ (or $r_2(t)$)
is the edge length opposite to angle t (or $x$) in the hyperbolic triangle with angles $(t, x, \pi-D(e))$.
Thus $\lim_{t\to 0}r_2(t)=\infty.$ Hence $\lim_{t\to 0}\frac d{dt}v(t,x)=0.$

And
$$\lim_{t\to 0}\frac d{dt}v(y_1-t,y_2)=\lim_{t\to 0} [y'_2\ln\tanh \frac{s_1(t)}2+  (y_1-t)'\ln\tanh
\frac{s_2(t)}2]=-\lim_{t\to 0}\ln\tanh \frac{s_2(t)}2,$$ where $s_1(t)$ (or $s_2(t)$)
is the edge length opposite to angle $y_1-t$ (or $y_2$) in the hyperbolic triangle with angles $(y_1-t,y_2,\pi-D(e'))$.
Thus $\lim_{t\to 0}s_2(t)$ are finite.
Hence $\lim_{t\to 0}\frac d{dt}v(y_1-t,y_2)$ is a positive finite number.

This is a contradiction.
Since we assume $\varphi_0$ is the maximum point, $\lim_{t\to 0}\frac d{dt}\mathcal{E}(\varphi_t)$
should be non-positive.

By the above argument we have show under the angle structure
$\varphi_0$, every triangle is a hyperbolic triangle. Hence
$\varphi_0\in\mathcal{A}.$ Therefore $\varphi_0$ is a critical
point of the energy function $\mathcal{E}.$ Since
$\sum_{e<f}\varphi(f,e)=\frac12\Phi(f),$ variables
$\{\varphi(f,e)\}$ are not free. By the method of Lagrange
multiplier, the angle structure $\{\varphi_0\}$ is a critical
point of the energy function if and only if there exists a
constant $c_f$ for each face $f$ such that $$0=\frac d {d
\varphi(f,e)}[\mathcal{E}(\varphi)-\sum_{f\in
F}c_f(\sum_{e<f}\varphi(f,e)-\frac12\Phi(f))]=\ln\tanh
\frac{r(f,e)}2-c_f.$$ This shows that the edge length $r(f,e)$
does not depend on any edges of the face $f$ but depends only on
the face $f$ itself. Thus we can take it to be the radius of the
face $f.$ This gives us a radius function $r:F\to(0,\infty)$,
equivalently a hyperbolic circle pattern.

By the argument above, a hyperbolic circle pattern induces an
angle structure which is a critical point of the energy function
$\mathcal{E}.$ Since the function $\mathcal{E}$ is strictly
concave down in the convex set $\mathcal{A}$, it has a unique
critical point which implies the uniqueness of the hyperbolic
circle pattern up to isometry.

\section*{Acknowledgement} The author would like to thank his advisor, Feng
Luo, for guidance and encouragement on this work.

\bibliographystyle{amsplain}

\end{document}